\documentclass[11pt]{article}

\usepackage{amsmath,amsbsy,amssymb,latexsym,amsthm}

\usepackage{natbib}

\newtheorem{Theorem}{Theorem}
\newtheorem{Definition}{Definition}
\newtheorem{Remark}{Remark}
\newtheorem{Lemma}{Lemma}
\newtheorem{Proposition}{Proposition}

\renewenvironment{abstract}{\begin{center}
\begin{minipage}[c]{12cm} {\begin{center}\bf Abstract \end{center}}} {\end{minipage}
\end{center}}
\newenvironment{keywords}{\begin{center}
\begin{minipage}[c]{12cm} {\bf Keywords:}} {\end{minipage}
\end{center}}
\newenvironment{msc}{\begin{center}
\begin{minipage}[c]{12cm} {\bf 2000 Mathematics Subject Classification:}} {\end{minipage}
\end{center}}

\begin{document}

\title{Isoperimetric problems on time scales\\
with nabla derivatives\thanks{This is a preprint
of an article whose final and definitive form will
be published in the \emph{Journal of Vibration and Control}. 
Date Accepted:  November 21, 2008.}}

\author{Ricardo Almeida\\
\texttt{ricardo.almeida@ua.pt}
\and
Delfim F. M. Torres\\
\texttt{delfim@ua.pt}}

\date{Department of Mathematics\\
University of Aveiro\\
3810-193 Aveiro, Portugal}

\maketitle


\begin{abstract}
We prove a necessary optimality condition
for isoperimetric problems under
nabla-differentiable curves.
As a consequence, the recent results
of [M.R. Caputo, A unified view of ostensibly disparate isoperimetric variational problems, Appl. Math. Lett. (2008),
doi:10.1016/j.aml.2008.04.004],
that put together seemingly dissimilar optimal control problems
in economics and physics, are extended to a generic time scale.
We end with an illustrative example of application
of our main result to a dynamic optimization
problem from economics.
\end{abstract}


\begin{msc}
49K30; 39A12.
\end{msc}


\section{Introduction}

The calculus on time scales is a recent field of mathematics, introduced by Bernd Aulbach and Stefan Hilger
[\citet{Hilger}], which unifies the theory of difference equations with the theory of differential equations. It has found applications in several fields that require simultaneous modeling of discrete and continuous data, in particular in control theory
[\citet{Barto:et:al}], [\citet{BartoPaw}], [\citet{MozBart}] and the calculus of variations [\citet{Atici}], [\citet{Bohner}],
[\citet{FerreiraTorres:2007}], [\citet{FerreiraTorres:a}],
[\citet{MalinowskaTorres}].

In Section~\ref{sec:prm} we present a short introduction to time scales and nabla derivatives. Section~\ref{sec:isoprob} is the main core of the paper: we prove
a necessary optimality condition for the isoperimetric problem on time scales with nabla-derivatives. Differently from [\citet{FerreiraTorres}], where the minimizing curve is assumed not to be an extremal of the constraint delta-integral,
here both normal (Theorem~\ref{normalcase}) and abnormal extremals (Theorem~\ref{thm:abn}) are considered. As a result of Theorem~\ref{thm:abn}, we extend the recent result of [\citet{Caputo}] to a generic time scale (Proposition~\ref{prop:ext:Caputo}). Finally,
in Section~\ref{sec:example} we illustrate
the application of Theorem~\ref{normalcase}
to an isoperimetric problem on time scales motivated by
[\citet{Atici}].


\section{Preliminaries}
\label{sec:prm}

For an introduction to time scales we refer the reader
to the comprehensive books [\citet{Bohner:Peterson}],
[\citet{Bohner-Peterson2}]. Here we just recall
the results and notation needed in the sequel.

By $\mathbb T$ we denote a time scale, \textrm{i.e.},
a nonempty closed subset of $\mathbb R$.
The \emph{backward jump operator}
$\rho:\mathbb{T}\rightarrow\mathbb{T}$ is defined by
$\rho(t)=\sup{\{s\in\mathbb{T}:s<t}\}$ if $t\neq \inf \mathbb{T}$,
and $\rho(\inf \mathbb{T})=\inf \mathbb{T}$. Analogously,
we define the \emph{forward jump operator}
$\sigma:\mathbb{T}\rightarrow\mathbb{T}$ by
$\sigma(t)=\inf{\{s\in\mathbb{T}:s>t}\}$ if $t\neq \sup \mathbb{T}$, and $\sigma(\sup \mathbb{T})=\sup
\mathbb{T}$.

A point $t\in\mathbb{T}$ is called \emph{right-dense},
\emph{right-scattered}, \emph{left-dense} and
\emph{left-scattered} if $\sigma(t)=t$, $\sigma(t)>t$, $\rho(t)=t$,
and $\rho(t)<t$, respectively. Let
$\mathbb{T}_\kappa$ be the set defined in the following way:
if $m$ is a right-scattered minimum of $\mathbb{T}$, then
$\mathbb{T}_\kappa=\mathbb{T}\setminus\{m\}$; if not,
$\mathbb{T}_\kappa= \mathbb{T}$.

A function $f:\mathbb{T}\rightarrow\mathbb{R}$ is called
\emph{nabla differentiable} at $t\in\mathbb{T}_\kappa$ if there
exists a number $f^{\nabla}(t)$ (called the \emph{nabla
derivative} of $f$ at $t$) such that for every $\epsilon>0$ there
exists some neighborhood $U$ of $t$ at $\mathbb{T}$ with
$$|f(\rho(t))-f(s)-f^{\nabla}(t)(\rho(t)-s)|\leq\epsilon|\rho(t)-s|$$
for all $s\in U$.
If $f$ is \emph{nabla differentiable} at all
$t\in\mathbb{T}_\kappa$, then we say that $f$
is \emph{nabla differentiable}.

When $\mathbb{T}=\mathbb{R}$, $f$ is nabla differentiable  at $t$
if and only if is differentiable at $t$. If
$\mathbb{T}=\mathbb{Z}$, then $f$ is always nabla differentiable
and
$$f^{\nabla}(t)=\frac{f(t)-f(\rho(t))}{t-\rho(t)}=f(t)-f(t-1).$$

Along the text we abbreviate $f \circ \rho$ by $f^\rho$ and
$[a,b]\cap \mathbb{T}$ by $[a,b]$. We use $\mathcal{C}^1([a,b],\mathbb R)$ to denote the set
$$
\mathcal{C}^1([a,b],\mathbb{R}) :=
\{ y:[a,b] \to \mathbb{R} \, | \, y^{\nabla} \text{ exists and is continuous on } [a,b]_{\kappa}\} \, .$$

We say that $F:\mathbb{T}\rightarrow\mathbb{R}$
is a \emph{nabla antiderivative} of
$f$ if
$$F^{\nabla}(t)=f(t), \qquad  \forall t \in \mathbb{T}_\kappa \, .$$ We then define the \emph{nabla integral}
of $f$ by
$$\int_{a}^{b}f(t)\nabla t:=F(b)-F(a) \, .$$

The following two formulas of nabla integration by parts hold:

$$\int_{a}^{b}f^\rho(t)g^{\nabla}(t)\nabla t=\left[(fg)(t)\right]_{t=a}^{t=b}-
\int_{a}^{b}f^{\nabla}(t)g(t)\nabla t$$
and
$$\int_{a}^{b}f(t)g^{\nabla}(t)\nabla t=\left[(fg)(t)\right]_{t=a}^{t=b}-
\int_{a}^{b}f^{\nabla}(t)g^\rho(t)\nabla t.$$

We end our brief review of the calculus on time scales
via nabla derivatives by recalling a fundamental lemma of the calculus of variations recently proved in
[\citet{MartinsTorres}]:

\begin{Lemma}
\label{fundamentallemma}
Let $f(\cdot) \in C([a,b], \mathbb{R})$. If $\int_{a}^{b}f(t)\eta^{\nabla}(t)\nabla t=0$
for every curve $\eta(\cdot)
\in \mathcal{C}^1([a,b], \mathbb{R})$ satisfying $\eta(a)=\eta(b)=0$, then $f(t)=c$,
$c \in \mathbb{R}$, for all $t\in [a,b]_\kappa$.
\end{Lemma}


\section{Main results}
\label{sec:isoprob}

We study the isoperimetric problem on time scales
with a nabla-integral constraint both for normal
and abnormal extremizers. The problem consists
of minimizing or maximizing
\begin{equation}
\label{eq:f:p:i}
I[y(\cdot)]=\int_{a}^{b}f(t,y^\rho(t),y^\nabla(t))\nabla t
\end{equation}
in the class of functions $y(\cdot) \in \mathcal{C}^1([a,b],\mathbb R)$
satisfying the boundary conditions
\begin{equation}
\label{eq:bc}
y(a)=\alpha \quad \mbox{ and } \quad  y(b)=\beta
\end{equation}
and the nabla-integral constraint
\begin{equation}
\label{eq:ic}
J[y(\cdot)]=\int_{a}^{b}g(t,y^\rho(t),y^\nabla(t))
\nabla t =\Lambda \, ,
\end{equation}
where $\alpha$, $\beta$, and $\Lambda$ are given real numbers.
We assume that functions
$(t,x,v) \rightarrow f(t,x,v)$ and $(t,x,v) \rightarrow g(t,x,v)$
possess continuous partial derivatives with respect to the second and third variables, and we denote them by $f_x$, $f_v$, $g_x$,
and $g_v$.

\begin{Definition}
We say that
$y(\cdot) \in \mathcal{C}^{1}([a,b],\mathbb{R})$
is a (weak) local minimizer (respectively local
maximizer) for the isoperimetric problem
\eqref{eq:f:p:i}--\eqref{eq:ic}
if there exists
$\delta >0$ such that
$I[y(\cdot)] \leq I[\hat{y}(\cdot)]$
(respectively $I[y(\cdot)] \geq I[\hat{y}(\cdot)]$)
for all $\hat{y}(\cdot) \in \mathcal{C}^{1}([a,b], \mathbb{R})$ satisfying the boundary conditions \eqref{eq:bc},
the isoperimetric constraint \eqref{eq:ic}, and
$||\hat{y}^\rho(\cdot) - y^\rho(\cdot)||
+ ||\hat{y}^\nabla(\cdot) - y^\nabla(\cdot)|| < \delta$,
where $||\hat{y}(\cdot)|| :=\sup_{t \in [a,b]_{\kappa}}|\hat{y}(t)|$.
\end{Definition}

\begin{Definition}
We say that $y(\cdot) \in \mathcal{C}^1([a,b],\mathbb R)$ is an \emph{extremal} for $J[\cdot]$ if
$$g_v\left(t,y^\rho(t),y^\nabla(t)\right)
-\int_a^t
g_x\left(\tau,y^\rho(\tau),y^\nabla(\tau)\right)\nabla\tau=c$$
for some constant $c$ and for all $t \in [a,b]_{\kappa}$.
An extremizer (\textrm{i.e.}, a local minimizer
or a local maximizer) to problem \eqref{eq:f:p:i}--\eqref{eq:ic} that is not an extremal for $J[\cdot]$ is said to be a normal extremizer; otherwise (\textrm{i.e.}, if it is an extremal for $J[\cdot]$), the extremizer is said to be \emph{abnormal}.
\end{Definition}

\begin{Theorem}[necessary optimality condition for normal extremizers of \eqref{eq:f:p:i}--\eqref{eq:ic}]
\label{normalcase}
Let $\mathbb{T}$ be a time scale, $a,b \in \mathbb{T}$ with
$a < b$, and $y(\cdot)\in \mathcal{C}^1([a,b],\mathbb R)$.
Suppose that $y(\cdot)$ gives a local minimum or a local maximum
to the functional $I[\cdot]$
subject to the boundary conditions $y(a)=\alpha$ and $y(b)=\beta$
and the integral constraint $J[y(\cdot)]=\Lambda$,
where $\alpha$, $\beta$, and $\Lambda$
are prescribed real values. If $y(\cdot)$ is not an extremal for $J[\cdot]$, then there exists a real $\lambda$ such that
$$F_v^\nabla(t,y^\rho(t),y^\nabla(t))-F_x(t,y^\rho(t),y^\nabla(t))=0$$
for all $t \in [a,b]_{\kappa}$, where $F = f-\lambda g$.
\end{Theorem}

\begin{proof}
Consider a variation of $y(\cdot)$, say
$\hat{y}(\cdot)=y(\cdot)+\epsilon_1\eta_1(\cdot)+\epsilon_2\eta_2(\cdot)$,
where $\eta_i(\cdot)$
is a curve in $\mathcal{C}^1([a,b],\mathbb R)$ satisfying $\eta_i(a)=\eta_i(b)=0$, $i=1$, $2$. Define the real functions
$$\hat I(\epsilon_1,\epsilon_2):=I[\hat{y}(\cdot)] \quad \mbox{ and }
 \quad \hat J(\epsilon_1,\epsilon_2):=J[\hat{y}(\cdot)]-\Lambda \, .$$
Integration by parts gives
$$\begin{array}{ll}
\displaystyle\left.\frac{\partial \hat J}{\partial \epsilon_2}\right|_{(0,0)} & =
          \displaystyle\int_{a}^{b} \left(\eta_2^\rho g_x\left(t,y^\rho,y^\nabla\right)+
          \eta_2^\nabla g_v\left(t,y^\rho,y^\nabla\right)\right)
          \nabla t \\
          &= \displaystyle\int_{a}^{b} \left[\eta_2^\rho
          \left(\int_a^t g_x\left(\tau,y^\rho,y^\nabla\right)
           \nabla\tau \right)^\nabla
          +\eta_2^\nabla g_v\left(t,y^\rho,y^\nabla\right)\right]
          \nabla t \\
          &= \displaystyle\int_{a}^{b}\left(   - \displaystyle\int_a^t g_x\left(\tau,y^\rho,y^\nabla\right)
 \nabla\tau + g_v\left(t,y^\rho,y^\nabla\right)\right)
\eta_2^\nabla \nabla t \\
\end{array}$$
since $\eta_2(a)=\eta_2(b)=0$. By Lemma~\ref{fundamentallemma},
there exists a curve $\eta_2(\cdot)$ such that $\left.\frac{\partial \hat
J}{\partial \epsilon_2}\right|_{(0,0)}\not=0$. Since $\hat
J(0,0)=0$, by the implicit function theorem there exists a
function $\epsilon_2(\cdot)$ defined in a neighborhood of zero,
such that $\hat J(\epsilon_1,\epsilon_2(\epsilon_1))=0$,
\textrm{i.e.}, we may choose a subset of the variation curves
$\hat{y}(\cdot)$ satisfying the isoperimetric constraint.

In sum, $(0,0)$ is an extremal of $\hat I$ subject to the
constraint $\hat J=0$ and $\textrm{grad}\hat{J}(0,0)\not=0$.
By the Lagrange multiplier rule, there exists
some $\lambda$ such that
$\textrm{grad}(\hat I-\lambda \hat J)=0$. Analogously,
$$\displaystyle\left.\frac{\partial \hat I}{\partial \epsilon_1}\right|_{(0,0)} =
\int_{a}^{b}\left(   - \int_a^t f_x(\tau,y^\rho,y^\nabla) \nabla\tau +
f_v(t,y^\rho,y^\nabla)\right)\eta_1^\nabla \nabla t $$
and
$$\displaystyle\left.\frac{\partial \hat J}{\partial \epsilon_1}\right|_{(0,0)} =
\int_{a}^{b}\left(   - \int_a^t g_x(\tau,y^\rho,y^\nabla) \nabla\tau +
g_v(t,y^\rho,y^\nabla)\right)\eta_1^\nabla \nabla t. $$
Therefore,
\begin{multline}
\label{eq1}
\int_{a}^{b}\left[   -
\int_a^t f_x(\tau,y^\rho,y^\nabla) \nabla\tau + f_v(t,y^\rho,y^\nabla) \right. \\
- \left.\lambda \left(  - \int_a^t g_x(\tau,y^\rho,y^\nabla) \nabla\tau + g_v(t,y^\rho,y^\nabla)
\right)\right]\eta_1^\nabla \nabla t=0 \, .
\end{multline}
Since (\ref{eq1}) holds for any curve $\eta_1(\cdot)$,
again by Lemma~\ref{fundamentallemma},
\begin{multline}
\label{eq2}
- \int_a^t f_x(\tau,y^\rho,y^\nabla) \nabla\tau + f_v(t,y^\rho,y^\nabla) \\
- \lambda \left(  - \int_a^t g_x(\tau,y^\rho,y^\nabla) \nabla\tau + g_v(t,y^\rho,y^\nabla) \right)=const \, .
\end{multline}
Applying the nabla derivative to both sides
of equation (\ref{eq2}), we get
$$
- f_x(t,y^\rho,y^\nabla) + f_v^\nabla(t,y^\rho,y^\nabla) -\lambda \left(-g_x(t,y^\rho,y^\nabla)+
g_v^\nabla(t,y^\rho,y^\nabla) \right)=0 \, ,
$$
\textrm{i.e.},
$$
F_v^\nabla(t,y^\rho,y^\nabla) - F_x(t,y^\rho,y^\nabla)=0 \, ,
$$
where $F=f-\lambda g$.
\end{proof}

\begin{Theorem}[necessary optimality condition for normal and abnormal extremizers of problem \eqref{eq:f:p:i}--\eqref{eq:ic}]
\label{thm:abn}
Let $\mathbb{T}$ be a time scale, $a,b \in \mathbb{T}$ with
$a < b$. If $y(\cdot)\in \mathcal{C}^1([a,b],\mathbb R)$
is a local minimizer or maximizer for $I[\cdot]$ subject to the boundary conditions $y(a)=\alpha$ and $y(b)=\beta$ and
to the integral constraint $J[y(\cdot)]=\Lambda$,
then there exist two constants $\lambda_0$ and $\lambda$,
not both zero, such that
$$K_v^\nabla\left(t,y^\rho(t),y^\nabla(t)\right)
-K_x\left(t,y^\rho(t),y^\nabla(t)\right)=0$$
for all $t \in [a,b]_{\kappa}$,
where $K=\lambda_0 f-\lambda g$.
\end{Theorem}

\begin{Remark}
If $y(\cdot)$ is a normal extremizer, then one can choose
$\lambda_0 = 1$ in Theorem~\ref{thm:abn}, thus obtaining
Theorem~\ref{normalcase}. For abnormal extremizers,
Theorem~\ref{thm:abn} holds with $\lambda_0 = 0$.
The condition $(\lambda_0,\lambda) \ne 0$ guarantees
that Theorem~\ref{thm:abn} is a useful necessary condition.
\end{Remark}

\begin{proof}
Following the proof of Theorem \ref{normalcase},
since $(0,0)$ is an extremal of $\hat I$
subject to the constraint $\hat J=0$,
the abnormal Lagrange multiplier rule
(\textrm{cf.}, \textrm{e.g.}, [\citet{vanBrunt}])
asserts the existence of reals
$\lambda_0$ and $\lambda$, not both zero, such that
$\textrm{grad}(\lambda_0\hat I-\lambda \hat J)=0$. Therefore,
$$\lambda_0 \displaystyle\left.\frac{\partial \hat I}{\partial \epsilon_1}\right|_{(0,0)}-
\lambda \displaystyle\left.\frac{\partial \hat J}{\partial \epsilon_1}\right|_{(0,0)}=0$$
\begin{multline*}
\Leftrightarrow \int_{a}^{b}\left[  \lambda_0\left( - \int_a^t f_x(\tau,y^\rho,y^\nabla) \nabla\tau +
f_v(t,y^\rho,y^\nabla)\right) \right. \\
- \left.\lambda \left(  - \int_a^t g_x(\tau,y^\rho,y^\nabla) \nabla\tau + g_v(t,y^\rho,y^\nabla)
\right)\right]\eta_1^\nabla \nabla t=0 \, .
\end{multline*}
From the arbitrariness of $\eta_1(\cdot)$, we conclude by Lemma~\ref{fundamentallemma} that
\begin{multline*}
\lambda_0\left( - \int_a^t f_x(\tau,y^\rho,y^\nabla) \nabla\tau + f_v(t,y^\rho,y^\nabla)\right) \\
-\lambda \left(  - \int_a^t g_x(\tau,y^\rho,y^\nabla) \nabla\tau + g_v(t,y^\rho,y^\nabla)\right)= \text{const} \, .
\end{multline*}
The desired condition follows by nabla differentiation.
\end{proof}

In the recent paper [\citet{Caputo}] the classical isoperimetric problem (\textrm{i.e.}, problem \eqref{eq:f:p:i}--\eqref{eq:ic} with $\mathbb{T} = \mathbb{R}$) was studied for a very particular functional \eqref{eq:f:p:i} and a very particular constraint \eqref{eq:ic} that often occurs in economics and physics.
For that particular class of problems,
the extremal curve is shown to be a constant. Here we remark
that the main result of [\citet{Caputo}] is still valid on a generic time scale $\mathbb{T}$:

\begin{Proposition}
\label{prop:ext:Caputo}
Let $\mathbb{T}$ be a time scale, $a,b \in \mathbb{T}$ with
$a < b$, $f:\mathbb R^2\to\mathbb R$ be a $C^2$ function
and $\xi$ a real parameter. Assume that
$(x,\xi) \rightarrow f(x,\xi)$ satisfies the conditions $f_x(y^\rho(t),\xi)\not=0$, in some interval $[c,d]\subseteq[a,b]$, for all $\xi$ and for all admissible function $y(\cdot)$, and $f_{xx}(y^\rho(t),\xi)\leq0$, for all $t\in[a,b]$ and for all $\xi$ over an open interval containing all the admissible values of $y^\rho(\cdot)$. Then, there exists a unique solution $y(\cdot)$
for the isoperimetric problem
\begin{equation}
\label{eq:PIC}
\begin{gathered}
I[y(\cdot)]=\int_{a}^{b}f(y^\rho(t),\xi)\nabla t
\longrightarrow \textrm{extr}\\
J[y(\cdot)]=\int_{a}^{b}y^\rho(t)\nabla t =\Lambda \\
y(a)=\alpha \, , \quad y(b)=\beta \, .
\end{gathered}
\end{equation}
Moreover, $y^\rho(t) = \text{constant}$, $t \in [a,b]_\kappa$.
\end{Proposition}

\begin{proof}
The proof follows from Theorem~\ref{thm:abn} and is the same,
\emph{mutatis mutandis}, than the one of Proposition~1 of [\citet{Caputo}].
\end{proof}

\begin{Remark}
Proposition~\ref{prop:ext:Caputo} asserts that if
$y(\cdot)$ is a solution of problem \eqref{eq:PIC}, then
$y(t)$ is constant in $[a,\rho(b)] \cap \mathbb T$, say $y(t)
= c$. If $\mathbb T=\mathbb R$, then by the isoperimetric constraint, one has $y^\rho(t)=y(t)=\Lambda/(b-a)$. The constant
$c = \Lambda/(b-a)$ is precisely the one obtained in [\citet{Caputo}].
\end{Remark}

\begin{Remark}
Let $\mathbb T=\mathbb Z$. Then, $y(t)$ is constant in $\{a,\ldots,b-1 \}$. By the isoperimetric constraint,
$$\int_{a}^{b}y^\rho(t) \nabla t = \sum_{t=a+1}^b y^\rho(t) =  \sum_{t=a}^{b-1} y(t) =(b-a)c = \Lambda \, .$$
Similar to $\mathbb{T} = \mathbb{R}$, one has $c=\Lambda/(b-a)$.
\end{Remark}


\section{An example}
\label{sec:example}

We exemplify the use of Theorem~\ref{normalcase} with an economic model. Assume we wish to maximize the functional
$$I[y(\cdot)]=\int_a^b u(y^\rho(t))\, l^t \nabla t,$$
where $y^\rho(t)$ is the consumption during period $t \in [a,b]$, $u(\cdot)$ is the one-period utility, $l \in (1/2,1)$ is the discount factor, and $I[\cdot]$ is the lifetime utility.
For a  detailed explanation of the economic model we refer the reader to [\citet{Atici}]. We consider the maximization problem subject to the constraint
$$\int_a^b t \, y^\nabla(t) \nabla t=\,const.$$
Let $g(t,x,v)=tv$. Since
$$g_v(t,y^\rho,y^\nabla)-\int_a^t g_x(\tau,y^\rho,y^\nabla)\nabla\tau=t,$$
there are no abnormal extremals for the problem.
The augmented Lagrange function is
$$F(t,x,v,\lambda)=u(x)l^t-\lambda tv.$$
By Theorem~\ref{normalcase}, $F_x=F_v^\nabla$, and so
$$u'(x)l^t=-(\lambda t)^\nabla=-\lambda.$$
It follows that $u(x)=\lambda l^{-t}/ \ln l +C$. If $u(\cdot)$
is invertible, the extremal is given by
$$y^\rho(t) = u^{-1}\left( \lambda l^{-t}/\ln l+C \right).$$


\section*{Acknowledgments}

The authors are grateful to the support
of the \emph{Control Theory Group} (cotg)
from the \emph{Centre for Research on Optimization and Control}
(CEOC) given by the \emph{Portuguese Foundation for Science and Technology} (FCT), cofinanced through the
\emph{European Community Fund} FEDER/POCI 2010.



\end{document}